\pdfoutput=1

\documentclass[reqno]{amsart}
\RequirePackage[l2tabu, orthodox]{nag} 

\usepackage{hyperref}

\usepackage{courier}
\usepackage{color}
\usepackage{amsmath,amsthm,amssymb}
\usepackage{cases} 
\usepackage[cal=boondox,scr=boondoxo]{mathalfa}

\usepackage{graphicx}
\usepackage[figurename=Fig.,labelsep=period]{caption}

\usepackage{textcomp} 
\usepackage{sistyle}
\usepackage{physics}
\usepackage{mathtools}
\mathtoolsset{showonlyrefs=false}

\usepackage{cite}

\theoremstyle{remark}

\newtheorem*{remark*}{Remark}

\begin{document}

\title[Reciprocal symmetry breaking in Pareto sampling]{Reciprocal symmetry breaking in\\ Pareto sampling}
\author[H.-S. Niwa]{H.-S. Niwa}
\date{}


\begin{abstract}
Let $W_1,\ldots,W_N$ be a sample of $\mathrm{Pareto}(\alpha)$ random variables normalized by their sum, such that $\sum_i W_i=1$.
The $W_i$ may represent
the weights of valleys in a spin glass (if $0<\alpha<1$),
or the frequency of different lineages (families) in a genealogy.
This paper considers a population in which there are $N$ individuals reproducing with $\mathrm{Pareto}(\alpha)$ offspring-number distribution ($1<\alpha<2$).
The probability of two randomly-chosen individuals being siblings,
$Y_2=\sum_i W_i^2$, gives the sample mean of the normalized size of families,
and its reciprocal gives the effective number of families (or reproducing lineages) in the population, $N_{\mathrm{e}}=1/Y_2$.
The typical sample mean is very different from the average over all possible samples, i.e. $Y_2$ is not a self-averaging quantity.
The typical $Y_2$ and its reciprocal do not vary with $N$ in opposite ways.
Non-self-averaging effects are crucial in understanding genetic diversity in mass spawning species such as marine fishes.
\end{abstract}

\maketitle

\section{Introduction}
Motivated by considering mass spawning species with type-III (exponential) survivorship curve,
many authors~\cite{SBB13,Arnason-Halldorsdottir2015,Niwa-etal2016}
have been studied patterns of genetic variation within marine populations
and obtained convincing results suggesting that
reproductive skew among individuals explains the pattern of coalescence of ancestral lineages.
While the recruitment process has exponential decay in survival probability,
the exponential amplification of the number of matured offspring (or siblings) in a family
compensates the exponentially small probability of their surviving to reproductive maturity.
The combination of these two exponentials leads to power laws in the offspring-number distribution~\cite{Reed-Hughes2002,Niwa-etal2017}.

Since annual recruitment
is calculated by summing random offspring numbers,
when the offspring-number distribution is broad,
the sum deserves serious consideration,
i.e. the system requires two different kinds of averages.
One might take the average of relevant properties (or variables) over the whole population existing at any given time.
However, these averages may fluctuate in time even for very large populations.
One can thus envisage taking
the time average of these population averages over a very long time stretch,
which, if some sort of ergodic property is assumed, may be represented by the average over all possible realizations of the stochastic reproduction process.
Such a stochastic nature of sums of a large number of random variables goes under the name of ``lack of self-averaging''
\cite{Derrida-PhysD97}.

This paper studies the distribution properties of family (or sibship) sizes
in the population, where their reproductive success follows a power-law distribution.
Let $W_i$ be the relative frequency or weight of the $i$-th family
(satisfying $\sum_i W_i=1$).
The weight $W_i$ is the probability for a given individual to be found in the $i$-th family,
so the sum of squared weights of families, $Y_2=\sum_i W_i^2$, is the expected weight of the family containing it.
The $Y_2$ gives the sample mean weight of families,
and its reciprocal $N_{\mathrm{e}}{\,}(=Y_2^{-1})$ gives the effective number of families (or reproducing lineages) in the population
\cite{Wright1931}.
The effective population size $N_{\mathrm{e}}$ is a measure of how many individuals contribute to the next generation.
The $Y_2$ also gives the probability of two randomly-chosen individuals being siblings (i.e. the coalescence probability).
Write $\rho(w)$ for the distribution of the weights, such that
$\rho(w)\dd{w}$ is defined as the average number of families with weights between $w$ and $w+\dd{w}$ among an infinite number of replicate populations each undergoing the same reproduction process.
Note that $w\rho(w)$ is the probability that a randomly chosen individual belongs to a family of weight $w$.
From the knowledge of $\rho(w)$,
I compute the probability distributions of the sample mean ($Y_2$) and its reciprocal ($N_{\mathrm{e}}$).

Closely related questions concerning the moments of the distribution $\rho(w)$ have been studied in genealogical, or coalescent, models.
Schweinsberg \cite{Schweinsberg2003} built coalescents from a population model with power-law offspring-number distribution.
Huillet \cite{Huillet2014} derived coalescents from sampling from a power-law distribution,
including size-biasing on the total recruitment effects.

The non-self-averaging effects are present in a large variety of disordered systems in statistical physics, in particular of spin glasses
\cite{MezardPSTV84,Derrida-PhysD97}.
Derrida and Peliti \cite{Derrida-Peliti91} computed the genealogy statistics under the Wright-Fisher model,
and by exploiting the equivalence with spin glasses they showed that similar non-self-averaging effects occur.
The random structure of family trees was analyzed to show that
the mean distance between individuals (number of generations from the common ancestor) fluctuates on time scale of the order of $N$ (population size) generations
\cite{Serva2004,Serva2005}.
Coalescent processes arise in a natural way from spin-glass models
\cite{Bolthausen-Sznitman98},
which allows one
to make related predictions concerning the non-self-averaging properties of the genealogies of evolving populations
\cite{Brunet-Derrida2009,Brunet-Derrida2011,Brunet-Derrida2013}.

In population genetics the non-self-averaging effects are observed,
when considering heterozygosity $H$ over many realizations of the neutral evolutionary process.
The standard deviation of $H$, calculated under the Wright-Fisher model,
is of the same order of magnitude as the mean for large populations with small mutation rates
\cite{Stewart76}.
The large variations in $H$ are to be expected between different gene loci,
which was observed in \textit{Drosophila melanogaster}
\cite{Higgs-PRE95}.
Tajima's estimator~\cite{Tajima83} of population-scaled mutation rate ($\propto N_{\mathrm{e}}$) has a similar non-self-averaging property, even when sampling infinitely many loci
\cite{King2018}.

In ecology, $H$ is the probability that two randomly selected individuals are of different species.
Based on the log-series distribution of the species abundances \cite{Fisher-etal43} in a (neutral) community,
the variance of $H$ does not go to zero even in the limit of large community size~\cite{Watterson74}.

Usually, demographic stochasticity has effects on small populations.
Variances of family-size frequencies in a large population average out, so that such microscopic fluctuations at the individual level may not be extracted from macroscopic (i.e. population-level) measurements (e.g. interannual recruitment variability).
This paper explains that,
as a consequence of non-self-averaging effects,
macroscopic measurements may give access to microscopic fluctuations (and vice versa), which is important for marine ecological applications.

This paper is organized as follows.
In \S\ref{section:population model}, after providing a population model,
I show that the weights of the families are not self-averaging.
In \S\ref{section:distribution-Y2},
the probability distributions of $Y_2$ and $N_{\mathrm{e}}$ are computed.
In \S\ref{section:Derrida-trick},
these probability distributions are obtained numerically.
In \S\ref{sect:application}, as an application example, I explain how genetic diversity varies with annual recruitment strength.
The analysis is based on the asymptotic (i.e. large-population limit) behavior of the model.

\section{Pareto sampling}\label{section:population model}
\subsection{Population model}
Consider a population with a fixed number $N$ of individuals reproducing asexually.
Each individual $i\,(=1,\ldots,N)$ in a given generation is assigned a random value $X_i$ of reproductive success.
The variables $X_1,\ldots,X_N$ are independent and identically distributed copies of $X$ with Pareto density
\begin{equation}\label{eqn:Pareto-density}
 f_X(x)=\alpha x^{-\alpha-1}
\end{equation}
and cumulative distribution function $F_X(x)=1-x^{-\alpha}$
(with $\alpha>0$ and $x\geq 1$).
Let $\mu$ be the mean; one has $\mu=\alpha/(\alpha-1)$ when $\alpha>1$, and $\mu=\infty$ when $\alpha\leq 1$.
Upon normalizing $X_1,\ldots,X_N$ by their sum
\begin{equation*}
 R_N=\sum_{i=1}^{N} X_i,
\end{equation*}
one defines the weight $W_i$ of the term $X_i$ in the sum as
\begin{equation*}
 W_i=X_i/R_N
\end{equation*}
($i=1,\ldots,N$).
Each $W_i$ gives the probability of reproductive success of individual~$i$.
Given the population at some generation, for each individual at the following generation, one chooses at random with probability $W_i$ one parent $i\in\{1,\ldots,N\}$.
Any generation is replaced by a new one.
The values $X_i$ are drawn afresh in each generation.
Note that when the $W_i$'s are identical, the sampling procedure is equivalent to Wright-Fisher sampling.

The $X_i$ gives an analog of the number of potential offspring (i.e. surviving young to reproductive maturity) of individual $i$,
and the sum $R_N$ corresponds to the annual recruitment, i.e. the total number of offspring entering the (potentially reproductive) population.
Then, a fixed number $N$ of reproducing individuals are chosen at random among $R_N$ individuals of that generation.
When $\alpha>1$, by the law of large numbers, one has $R_N>N$.

In this paper I mainly consider the case $1<\alpha<2$.

\subsection{Domination by the largest term}
It is well known that
$(R_N-\mathrm{E}[R_N])/N^{1/\alpha}$ with $1<\alpha<2$ has a maximally asymmetric $\alpha$-stable distribution for large $N$
\cite{Levy1937},
where $\mathrm{E}[R_N]=\mu N$.
The width of the distribution of the sum $R_N$ (i.e. the typical value of the difference $R_N-\mathrm{E}[R_N]$) is of order $N^{1/\alpha}$,
while the variance $\mathrm{E}[R_N^2]-\mathrm{E}[R_N]^2$ is infinite.
$\mathrm{E}[\,\cdot\,]$ denotes an average over all possible samples (realizations).

Define $X_{1,N}\geq X_{2,N}\geq\cdots\geq X_{N,N}$ by ranking in decreasing order the values encountered among the $N$ terms of the sum $R_N$.
When $1<\alpha<2$, one has
\cite{ZKS05}
\begin{align*}
 \mathrm{E}\qty[X_{1,N}]&=\frac{N!{\,}\mathrm{\Gamma}(1-1/\alpha)}{\mathrm{\Gamma}(N+1-1/\alpha)}
 =\mathrm{\Gamma}(1-1/\alpha)N^{1/\alpha}\\
 \mathrm{E}\qty[X_{2,N}]&=\frac{\alpha-1}{\alpha}\mathrm{E}\qty[X_{1,N}]
\end{align*}
and while $X_{1,N}$ has an infinite second moment, one has
\begin{align*}
 \mathrm{E}\qty[X_{2,N}^2]&=\frac{N!{\,}\mathrm{\Gamma}(2-2/\alpha)}{\mathrm{\Gamma}(N+1-2/\alpha)}
 =\mathrm{\Gamma}(2-2/\alpha)N^{2/\alpha}\\
 \mathrm{E}\qty[X_{1,N}X_{2,N}]&=\frac{\alpha}{\alpha-1}\mathrm{E}\qty[X_{2,N}^2]
\end{align*}
for large $N$.
The rescaled random variable
$X_{1,N}/N^{1/\alpha}$
also has a $\mathrm{Pareto}(\alpha)$ distribution
and $X_{1,N}\gtrsim N^{1/\alpha}$
\cite{Bouchaud-Georges90,Hofstad2016}.
Importantly, all but the largest order statistics have finite second moment.
Therefore, the sum $R_{2,N}{\,}(=\sum_{i=2}^N X_{i,N})$ of the $(N-1)$ lower order statistics converges
to a normally distributed random variable with first two moments given by
\begin{align*}
 \mathrm{E}\qty[R_{2,N}] &=
 \frac{\alpha}{\alpha-1}\qty(N-\mathrm{\Gamma}(2-1/\alpha)N^{1/\alpha})
 =\mathrm{E}[R_N]\\
 \mathrm{E}\qty[R_{2,N}^2] &=
  \mathrm{E}\qty[\sum_{i=2}^N X_{i,N}^2]
  +\mathrm{E}\qty[\sum_{i\neq j} X_{i,N}X_{j,N}]
  \sim N^{2/\alpha}
\end{align*}
where
\begin{equation*}
 \mathrm{E}\qty[\sum_{i=2}^N X_{i,N}^2]=
  \frac{\alpha}{\alpha-2}
  \qty(N-\mathrm{\Gamma}(2-2/\alpha)N^{2/\alpha})
  =\frac{\alpha}{2-\alpha}\mathrm{E}\qty[X_{2,N}^2],
\end{equation*}
and
\begin{align*}
 \mathrm{E}\qty[X_{2,N}X_{3,N}] &=
  \frac{2\alpha}{2\alpha-1}\mathrm{E}\qty[X_{3,N}^2]\\
 \mathrm{E}\qty[X_{2,N}X_{4,N}] &=
  \frac{6\alpha^2}{(2\alpha-1)(3\alpha-1)}\mathrm{E}\qty[X_{4,N}^2]
\end{align*}
etc.
Accordingly, the statistical variation of the sum $R_N$ is dominated by its largest term $X_{1,N}$,
so the fraction $(R_N-\mathrm{E}[R_N])/R_N$ can be linked to one parent.
The concept of the statistical domination by the largest term is especially useful when describing processes with large deviations, as I show later.

\subsection{Moments of weights}
The fluctuations of the weights $W_i$ of the term $X_i$ in the sum $R_N$ can be described by considering their moments.
When the sum of the $k$-th power of weights ($k\geq 0$),
\begin{equation*}
 Y_k = \sum_{i=1}^N W_i^k,
\end{equation*}
is averaged over the $X_i$'s, one gets the moments $\mathrm{E}\qty[Y_k]$.
Obviously one has $\mathrm{E}[Y_0]=N$ and $\mathrm{E}[Y_1]=1$.
Using the following identity
\begin{equation*}
 \frac{\sum_{i=1}^N X_i^k}{\qty(\sum_{j=1}^N X_j)^k}=
 \int_0^{\infty}
 \frac{s^{k-1}\dd{s}e^{-s\sum_{j=1}^N X_j}}{\mathrm{\Gamma}(k)}
 \sum_{i=1}^N X_i^k
\end{equation*}
(this is a direct consequence of the Euler integral for the gamma function),
one can calculate $\mathrm{E}\qty[Y_k]$ in the large-$N$ limit,
\begin{equation}\label{eqn:moment-Y}
 \mathrm{E}\qty[Y_k]
 =
 \frac{N}{\mathrm{\Gamma}(k)}\int_0^{\infty} s^{k-1}\dd{s}
 \mathrm{E}\qty[e^{-sX_1}]^{N-1} \mathrm{E}\qty[X_2^k e^{-sX_2}],
\end{equation}
where
the integral is dominated by the small $s$ behavior.
I refer to~\cite{Derrida-PhysD97,Huillet2014}.
One sees, via integration by parts $\lfloor\alpha\rfloor+1$ times,
that for small $s$
\begin{equation}\label{eqn:exp-tx}
 \mathrm{E}\qty[e^{-sX}]=\sum_{0\leq k<\alpha}\frac{(-s)^{k}}{{k}!}\mathrm{E}\qty[X^{k}]-s^{\alpha}\mathrm{\Gamma}(1-\alpha)
\end{equation}
for $\alpha\neq 1,2,\ldots$,
where $\lfloor\cdot\rfloor$ denotes the integer part of the argument.
It is easy to check that
\begin{align*}
 & \mathrm{E}\qty[e^{-sX}] = \int_{1}^{\infty} e^{-sx}f_X(x)\dd{x}\\
 & \mbox{ }=
 \begin{dcases}
  1-\alpha\int_0^{\infty} x^{-1-\alpha}\qty(1-e^{-sx})\dd{x}&
  (0<\alpha<1)\\
  e^{-s}+s\int_{1}^{\infty}\qty(F_X(x)-1)\dd{x}
  -s\int_{1}^{\infty}\qty(F_X(x)-1)\qty(1-e^{-sx})\dd{x}&
  (\alpha>1)
 \end{dcases}
\end{align*}
where the third term on the last line reduces to
\begin{equation*}
 s\int_0^{\infty} x^{-\alpha}\qty(1-e^{-sx})\dd{x}
  =
  \frac{s^2}{\alpha-1}\int_0^{\infty} x^{-\alpha+1} e^{-sx}\dd{x}
\end{equation*}
for $1<\alpha<2$, and to
\begin{equation*}
 \frac{s\qty(1-e^{-s})}{\alpha-1}+
  \frac{s^2 e^{-s}}{(\alpha-1)(\alpha-2)}
  -\frac{s^3}{(\alpha-1)(\alpha-2)} \int_{1}^{\infty} x^{-\alpha+2} e^{-sx}\dd{x}
\end{equation*}
for $\alpha>2$.
Take note that, for small $s$, only large values of $x$ contribute to the integral.
One views the boundary case $\alpha=1$ (resp. $\alpha=2$, etc.) as the limiting critical case of Eq.\eqref{eqn:exp-tx} as
$\alpha\to 1+0$ (resp. $\alpha\to 2+0$, etc.).
So one gets
\begin{equation*}
 \mathrm{E}\qty[e^{-sX}] = 1-s+s\ln s
\end{equation*}
for $\alpha=1$, and
\begin{equation*}
 \mathrm{E}\qty[e^{-sX}] = 1-s\mu+\frac{s^2}{2}-s^2\ln s
\end{equation*}
for $\alpha=2$, etc.
Thus, one can see that for $k>\alpha$
\[
 \mathrm{E}\qty[X^k e^{-sX}]
 = \alpha s^{\alpha-k}\mathrm{\Gamma}(k-\alpha).
\]

When $0<\alpha<1$, Eq.\eqref{eqn:moment-Y} gives the equation~(11) of \cite{Derrida-PhysD97}.

When $1\leq\alpha<2$, one obtains for $k\geq 2$
\begin{equation*}
 \mathrm{E}\qty[Y_k] =
 c_{\scalebox{0.55}{$N$}}
 \frac{\mathrm{Beta}(k-\alpha,\alpha)}{\mathrm{Beta}(2-\alpha,\alpha)}
\end{equation*}
with scaling constant
\begin{equation*}
 c_{\scalebox{0.55}{$N$}} =
  \begin{dcases}
   \frac{\alpha\mathrm{Beta}(2-\alpha,\alpha)}{\mu^{\alpha}N^{\alpha-1}} & (1<\alpha<2)\\
    \qty(\ln N)^{-1}& (\alpha=1)
   \end{dcases}
\end{equation*}
where $\mathrm{Beta}(a,b)=\mathrm{\Gamma}(a)\mathrm{\Gamma}(b)/\mathrm{\Gamma}(a+b)$ is the beta function.

When $\alpha\geq 2$, one gets
\begin{numcases}{\mathrm{E}\qty[Y_k] =}
 \frac{\alpha\mathrm{Beta}(k-\alpha,\alpha)}{\mu^{\alpha}N^{\alpha-1}} & ($\alpha<k$)\nonumber\\
 \label{eqn:ell-moment-Kingman}
 \frac{\mathrm{E}[X^k]}{\mu^k N^{k-1}} & ($2\leq k<\alpha$)
\end{numcases}
and
\begin{equation}\label{eqn:int-alpha-moment-Kingman}
\mathrm{E}\qty[Y_{\alpha}]=\frac{\alpha\ln N}{\mu^{\alpha}N^{\alpha-1}}
\end{equation}
for integer $\alpha=2,3,\ldots$.

The correlations between the $Y_k$'s can also be calculated as
\begin{align*}
 \mathrm{E}\qty[Y_k Y_{\ell}] &=
 \mathrm{E}\qty[\sum_{i=1}^N W_i^{k+\ell}]
 +\mathrm{E}\qty[\sum_{i\neq j} W_i^k W_j^{\ell}]\\
 &= \frac{N}{\mathrm{\Gamma}(k+\ell)}\int_0^{\infty}
 s^{k+\ell-1}\dd{s}
 \mathrm{E}\qty[e^{-sX_1}]^{N-1} \mathrm{E}\qty[X_2^{k+\ell} e^{-sX_2}]\\
 &{\quad}+\frac{N^2}{\mathrm{\Gamma}(k+\ell)}\int_0^{\infty}
 s^{k+\ell-1}\dd{s}\mathrm{E}\qty[e^{-sX_1}]^{N-2}
 \mathrm{E}\qty[X_2^k e^{-sX_2}]\mathrm{E}\qty[X_3^{\ell} e^{-sX_3}].
\end{align*}
When $1\leq\alpha<2$, one obtains the $(k,\ell)$-th moment ($k,\ell\geq 2$) of weights of two different families,
\begin{equation}\label{eqn:2-simul-merger}
 \mathrm{E}\qty[\sum_{i\neq j} W_i^k W_j^{\ell}]
  =
  \mathrm{E}\qty[Y_k]\mathrm{E}\qty[Y_{\ell}]\frac{\mathrm{Beta}(k,\ell)}{\mathrm{Beta}(\alpha,\alpha)}.
\end{equation}
Therefore,
the correlations
$\mathrm{E}\qty[Y_k Y_{\ell}]$
dominate over the factorized terms
$\mathrm{E}\qty[Y_k]\mathrm{E}\qty[Y_{\ell}]$.

When $0<\alpha<1$,
the correlations $\mathrm{E}\qty[Y_k Y_{\ell}]$, as well as the moments $\mathrm{E}\qty[Y_k]$, depend only on $\alpha$, $k$ and $\ell$, and become independent of $N$;
see the equation~(12) of \cite{Derrida-PhysD97}.

\subsection{Fluctuation dominance}
To characterize fluctuations in the weights of families,
consider the relative fluctuations (i.e. coefficient of variation) of the $Y_2$,
\begin{equation*}
 \mathrm{CV}[Y_2]=
  \sqrt{\mathrm{E}[(Y_2-\mathrm{E}[Y_2])^2]}/\mathrm{E}[Y_2].
\end{equation*}
When $1<\alpha<2$, one has
\begin{equation*}
 \mathrm{CV}[Y_2]=
  c_{\scalebox{0.55}{$N$}}^{-1/2}\sqrt{(3-\alpha)(2-\alpha)/6},
\end{equation*}
which diverges in $\order{N^{(\alpha-1)/2}}$.
There is no approximately deterministic property of the $Y_2$-distribution at large population sizes,
as the fluctuations of the $Y_2$ dominate the average value.
Each realization of the $Y_2$ may be very different from its other realizations.
The fact that
$\mathrm{E}\qty[Y_k^2]/\mathrm{E}\qty[Y_k]^2=\order{N^{\alpha-1}}$
indicates that the $Y_k$'s for $k\geq 2$ are non-self-averaging.

When $\alpha\geq 2$,
one sees that $\mathrm{CV}[Y_2]=\order{N^{(3-\alpha)/2}}$ for $2<\alpha<4$, or $\mathrm{CV}[Y_2]=\order{\sqrt{N}/\ln N}$ for $\alpha=2$.
Although one ignores the probability that three or more randomly chosen individuals are siblings, i.e.
$\lim_{N\to\infty}\mathrm{E}[Y_k]/\mathrm{E}[Y_2]=0$ for all $k\geq 3$,
one cannot ignore the relative fluctuations in the probability of being siblings
when $2\leq\alpha\leq 3$.
Then, the probability of having four siblings from a family is greater than or equal to the probability of having two pairs of siblings from two different families.

\subsection{Distribution of the weights of families}
Since all the moments of weights are known,
the function $\mathrm{E}[Y_k]^{-1} w^k \rho(w)$ with $k\geq 2$ is known
\cite{Hausdorff1923},
which is the $\mathrm{Beta}(k-\alpha,\alpha)$ distribution
on $0\leq w\leq 1$.
For $1\leq\alpha<2$, the distribution
\begin{equation}\label{eqn:freq-spectrum-beta}
 \rho(w) =
  c_{\scalebox{0.55}{$N$}}
  \frac{w^{-\alpha-1}(1-w)^{\alpha-1}}{\mathrm{Beta}(2-\alpha,\alpha)}
\end{equation}
diverges like $w^{-\alpha-1}$ for small $w>0$.
Although there are a large number of very small families, the expected number of siblings of an individual,
$Nc_{\scalebox{0.55}{$N$}}$,
grows as $N^{2-\alpha}$ for large $N$.
The functions $\rho(w)$ and $w\rho(w)$ are not integrable
on $0\leq w\leq 1$.
If introducing a cut-off
\begin{equation}\label{eqn:cut-off}
 \varepsilon_{\scalebox{0.55}{$N$}}=(\mu N)^{-1}
\end{equation}
in the region of small $w$,
the total number of families having a weight larger than $\varepsilon_{\scalebox{0.55}{$N$}}$ is
\begin{equation*}
 \int_{\varepsilon_{\scalebox{0.45}{$N$}}}^1\rho(w) \dd{w}
 =
 \frac{c_{\scalebox{0.55}{$N$}}{\,}\varepsilon_{\scalebox{0.55}{$N$}}^{-\alpha}}{\alpha\mathrm{Beta}(2-\alpha,\alpha)} = N,
\end{equation*}
and their total weight is
\begin{equation*}
 \int_{\varepsilon_{\scalebox{0.45}{$N$}}}^1 w\rho(w) \dd{w}
  =
  \frac{c_{\scalebox{0.55}{$N$}}{\,}\varepsilon_{\scalebox{0.55}{$N$}}^{1-\alpha}}{(\alpha-1)\mathrm{Beta}(2-\alpha,\alpha)} = 1,
\end{equation*}
as $N\to\infty$.

Let $\rho(w,w')$ be the joint distribution (or the correlation function) of weights of two families,
\begin{equation*}
 \rho(w,w') = \rho(w)\delta\qty(w-w') + \rho^{\ast}(w,w').
\end{equation*}
The second term captures the average number of pairs of different families having weights $w$ and $w'$,
\begin{equation*}
 \rho^{\ast}(w,w') =
 c_{\scalebox{0.55}{$N$}}^2
 \frac{(ww')^{-\alpha-1} (1-w-w')^{2\alpha-1}}{\mathrm{Beta}(2-\alpha,\alpha)^2}{\,}\Theta\qty(1-w-w'),
\end{equation*}
which is extracted from Eq.\eqref{eqn:2-simul-merger} for $1\leq\alpha<2$,
where $\Theta(\cdot)$ is the Heaviside step function.
I refer to~\cite{MezardPSTV84,Derrida-Toulouse85}.
The correlated probability of finding two families of weights $w$ and $w'$ is given by $ww'\rho(w,w')$.

For $\alpha\geq 2$,
while when $2\leq k\leq\alpha$, the functions $w^k\rho(w)$ with
\begin{equation}\label{eqn:freq-spectrum-Kingman}
 \rho(w) =
  \frac{\alpha w^{-\alpha-1}(1-w)^{\alpha-1}}{\mu^{\alpha} N^{\alpha-1}}
\end{equation}
is not integrable
on $0\leq w\leq 1$,
one sees that
\begin{equation*}
 \int_{\varepsilon_{\scalebox{0.45}{$N$}}}^1 w^k\rho(w)\dd{w}
  =\mathrm{E}\qty[Y_k]
\end{equation*}
with $\mathrm{E}[Y_k]$ as in Eq.\eqref{eqn:ell-moment-Kingman} or~\eqref{eqn:int-alpha-moment-Kingman}.
Further, the function $\mathrm{E}[Y_2]^{-1}w^2\rho(w)$ reduces to the Dirac delta function,
as one sees as follows.
Define a function
\begin{equation*}
 \delta_N(w)=\mathrm{E}[Y_2]^{-1} w^2\rho(w)\Theta(w-\varepsilon_{\scalebox{0.45}{$N$}})
\end{equation*}
with $\rho(w)$ as in Eq.\eqref{eqn:freq-spectrum-Kingman}.
When $N\to\infty$, the $\delta_N(w)$ approximates the Dirac delta distribution at zero.

The joint distribution of the sequence $\qty{W_1,\ldots,W_N}$, in the case $0<\alpha<1$, has been studied in mathematics and physics.
The sequence of $W_1,\ldots,W_N$ in decreasing order, as $N\to\infty$, has the two-parameter Poisson-Dirichlet distribution with parameters $(\alpha,0)$
\cite{Pitman-Yor97}.
The distribution $\rho(w)$ and the joint distributions $\rho(w,w'), \rho(w,w',w''),\ldots$ of weights of two or more families are derived in \cite{MezardPSTV84}.
Although families are mostly concentrated around $w=0$,
these $w\simeq 0$ families do not contribute to the total weight, as
\begin{equation*}
 \int_0^w w'\rho(w')\dd{w'}\sim w^{1-\alpha}
\end{equation*}
with $\rho(w)$ given by the equation~(34) of \cite{MezardPSTV84}.
So any one of them has an extremely small weight.
The expected number of siblings of an individual is $(1-\alpha)N$ and
there are a few families with weights of $\order{1}$ in $1/N$.

\section{Distributions of
\texorpdfstring{$Y_2$}{TEXT}
and
\texorpdfstring{$N_{\mathrm{e}}$}{TEXT}
}\label{section:distribution-Y2}
Given the weights of families, $W_1,\ldots,W_N$,
write $W_{1.N}$ and $W_{2,N}$ for the weights of the largest and second largest families.
Let $\mathrm{\Pi}_{W_{1,N}}(w)$ (resp. $\mathrm{\Pi}_{W_{2,N}}(w)$) be the probability of the largest family (resp. the second largest family) having a weight $w$.
These probability distributions can be computed from the knowledge of the distribution $\rho(w)$ and the joint distributions of weights;
refer to~\cite{Derrida-Flyvbjerg87,Derrida-Flyvbjerg87-PhysA}.
If there is a family with weight $w>1/2$, this weight must be the largest one, and thus
$\mathrm{\Pi}_{W_{1,N}}(w)=\rho(w)$ for $w>1/2$.
Letting $\rho^{\ast}(w_1,\ldots,w_n)$ be a joint distribution of $n{\;(\geq 2)}$ weights of different families, and denoting
\begin{equation*}
 I_n(w)=\int_w^1\dd{v_1}\int_w^{v_1}\dd{v_2}\cdots\int_w^{v_{n-2}}\dd{v_{n-1}}\rho^{\ast}(v_1,v_2,\ldots,v_{n-1},w),
\end{equation*}
one has,
in the interval $1/(n+1)<w<1/n$,
\begin{equation*}
 \mathrm{\Pi}_{W_{1,N}}(w)=\rho(w)-I_2(w)+\cdots+(-1)^{n-1}I_n(w)
\end{equation*}
and
\begin{equation*}
 \mathrm{\Pi}_{W_{2,N}}(w)=\int_w^1\dd{v}\rho^{\ast}(v,w)-2I_3(w)+\cdots+(-1)^{n-2}(n-1)I_n(w).
\end{equation*}
The successive terms in these summations rapidly decrease for large $N$,
because $I_n(w)$ is an integral over more and more variables and diminishes in magnitude with increasing $n$.

Consider the case $1<\alpha<2$.
Since $W_{1,N}\simeq X_{1,N}/(\mu N)$, one has $W_{1,N}\gtrsim N^{1/\alpha-1}$.
Accordingly, one has, up to the leading term,
\begin{equation}\label{eqn:largest-W}
 \mathrm{\Pi}_{W_{1,N}}(w) = \rho(w)
\end{equation}
on $N^{1/\alpha-1}\lesssim w\leq 1$, and
\begin{equation}\label{eqn:2nd-largest-W}
 \mathrm{\Pi}_{W_{2,N}}(w) = \int_w^1\dd{v}\rho^{\ast}(v,w)
  =
  \frac{c_{\scalebox{0.55}{$N$}}^2{\,} w^{-2\alpha-1}(1-w)^{2\alpha-1}}{\alpha\mathrm{Beta}(2-\alpha,\alpha)}
\end{equation}
on $N^{1/\alpha-1}\lesssim w\leq 1$.
Integrating Eq.\eqref{eqn:largest-W}, one gets
\begin{equation*}
 \int_{N^{1/\alpha-1}/\mu}^1 \mathrm{\Pi}_{W_{1,N}}(w)\dd{w} = 1.
\end{equation*}
One then sees that the the second moment of the largest weight
\begin{equation*}
 \int_{N^{1/\alpha-1}/\mu}^1 w^2\mathrm{\Pi}_{W_{1,N}}(w)\dd{w}= c_{\scalebox{0.55}{$N$}}
\end{equation*}
for large $N$, which demonstrates that the sum of the squared weights of families, $Y_2=\sum_{i=1}^N W_i^2$, is dominated by its largest term $W_{1,N}$.

Now compute the probability of the random variable $Y_2$ taking a value $y$.
The probability distribution $\mathrm{\Pi}_{Y_2}(y)$ of $Y_2$ is defined as
\begin{equation*}
 \mathrm{\Pi}_{Y_2}(y)=\mathrm{E}\qty[\delta\qty(\sum_{i=1}^N W_i^2-y)].
\end{equation*}
Since the sums which contribute to the average are those in which one term dominates,
one may replace the sum with its maximal summand~\cite{MezardPSTV84}
\begin{equation*}
 \mathrm{\Pi}_{Y_2}(y) = \mathrm{E}\qty[\delta\qty(W_{1,N}^2-y)].
\end{equation*}
From the probability density $\mathrm{\Pi}_{W_{1,N}}(w)$ of the largest weight,
one obtains
\begin{equation}\label{eqn:distr-Y2}
 \mathrm{\Pi}_{Y_2}(y) =
  \int_{N^{1/\alpha-1}/\mu}^1\dd{w}\mathrm{\Pi}_{W_{1,N}}(w)\delta\qty(w^2-y)
  =\frac{c_{\scalebox{0.55}{$N$}} y^{-\alpha/2-1}\qty(1-\sqrt{y})^{\alpha-1}}{2\mathrm{Beta}(2-\alpha,\alpha)}
\end{equation}
on $N^{2(1/\alpha-1)}\lesssim y\leq 1$.
Integrating the right side of Eq.\eqref{eqn:distr-Y2} gives
\begin{equation*}
 \int_{N^{2(1/\alpha-1)}/\mu^{2}}^1 \mathrm{\Pi}_{Y_2}(y)\dd{y}=1,
\end{equation*}
and one sees that
\begin{equation*}
 \int_{N^{2(1/\alpha-1)}/\mu^{2}}^1 y{\,}\mathrm{\Pi}_{Y_2}(y)\dd{y}
  = c_{\scalebox{0.55}{$N$}}.
\end{equation*}
The most probable value of $Y_2$ is close to $N^{2(1/\alpha-1)}$, and therefore differs from its mean value $c_{\scalebox{0.55}{$N$}}\,(\sim N^{1-\alpha})$.

The same kind of argument applies to the probability distribution $\mathrm{\Pi}_{N_{\mathrm{e}}}(y)$ of $N_{\mathrm{e}}=Y_2^{-1}$, yielding
\begin{equation}\label{eqn:distr-Ne}
 \mathrm{\Pi}_{N_{\mathrm{e}}}(y)=
  \frac{c_{\scalebox{0.55}{$N$}} y^{\alpha/2-1}\qty(1-1/\sqrt{y})^{\alpha-1}}{2\mathrm{Beta}(2-\alpha,\alpha)}
\end{equation}
for $1\leq y\lesssim N^{2(1-1/\alpha)}$.
Integrating the right side of Eq.\eqref{eqn:distr-Ne} gives
\begin{equation*}
 \int_1^{\mu^2 N^{2(1-1/\alpha)}}\mathrm{\Pi}_{N_{\mathrm{e}}}(y)\dd{y}=1.
\end{equation*}
The function $\mathrm{\Pi}_{N_{\mathrm{e}}}(y)$ must break at $y\sim N^{2(1-1/\alpha)}$ and will decay rapidly when $y\to\infty$.
Moreover,
the most probable value of $N_{\mathrm{e}}$ is $(2-\alpha)^{-2}$ and
very different from the harmonic mean of $N_{\mathrm{e}}$'s
($c_{\scalebox{0.55}{$N$}}^{-1}\sim N^{\alpha-1}$)
over replicate populations.

\begin{remark*}
While one might naively expect an increase of the typical $N_{\mathrm{e}}$ with $N$,
the typical $N_{\mathrm{e}}$ shows no increase with population size $N$.
The paper presents the striking scaling behavior of the sample mean reciprocal.
The $Y_2$-distribution is of reciprocal symmetry breaking,
in the sense that
the typical sample mean and its reciprocal do not vary with population size in opposite ways.
The typical $Y_2$ decreases with the population size like $N^{2(1/\alpha-1)}$,
while the typical reciprocal of $Y_2$ is independent of $N$.
Such counter-intuitive behavior emerges in some broad distributions
\cite{Romeo2003}.
\end{remark*}

\section{Numerical reconstruction of probability distributions}\label{section:Derrida-trick}
\begin{figure}[hbt!]
 \centering
 \begin{tabular}{ll}
  \small{(a)} & \small{(b)}\\
  \includegraphics[height=.3\textwidth,viewport=0 0 360 232]{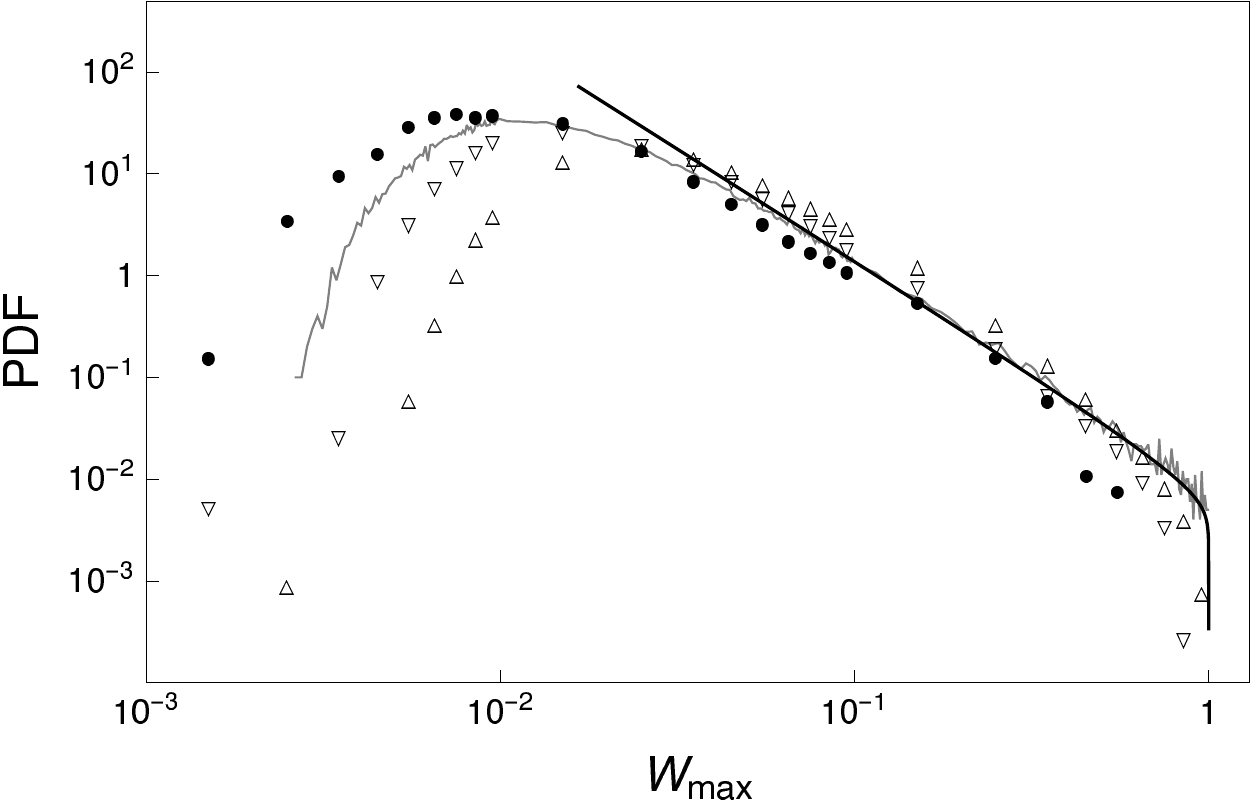}&
      \includegraphics[height=.3\textwidth,viewport=0 0 360 232]{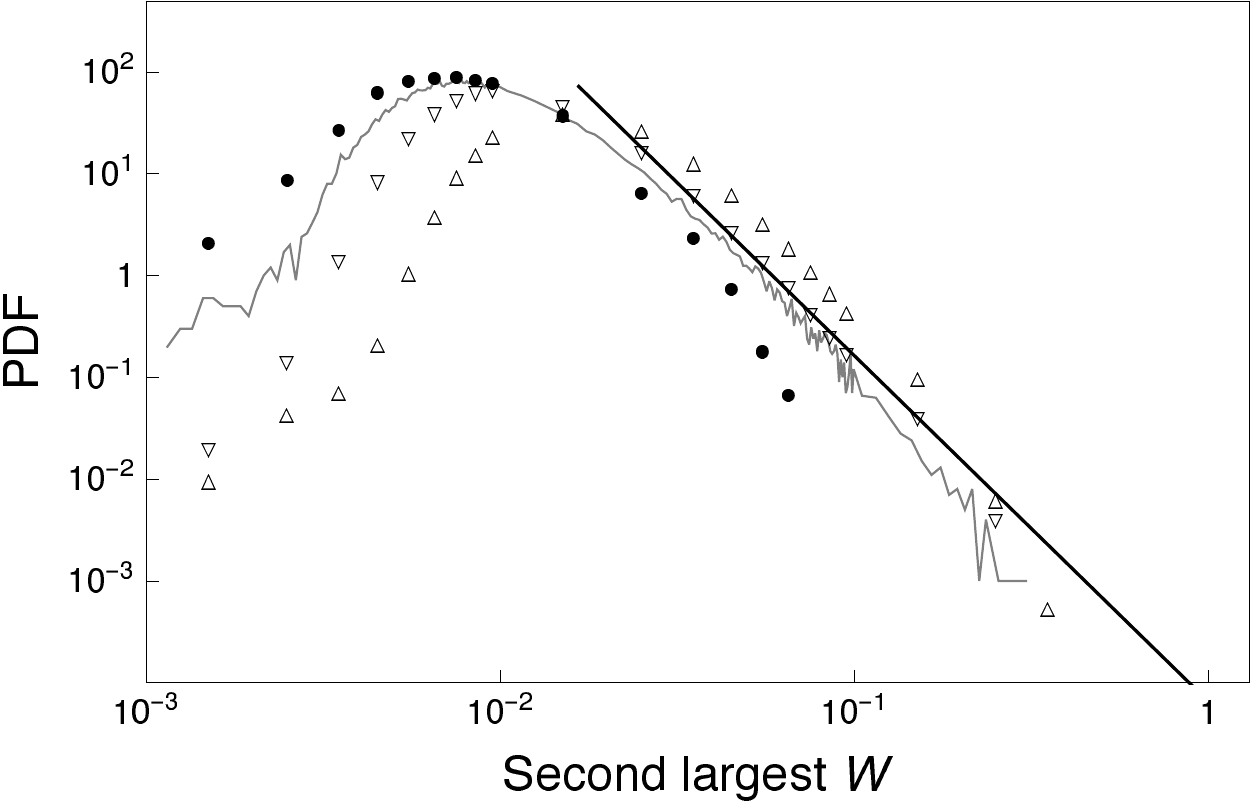}\\
  \small{(c)} & \small{(d)}\\
  \includegraphics[height=.3\textwidth,viewport=0 0 360 231]{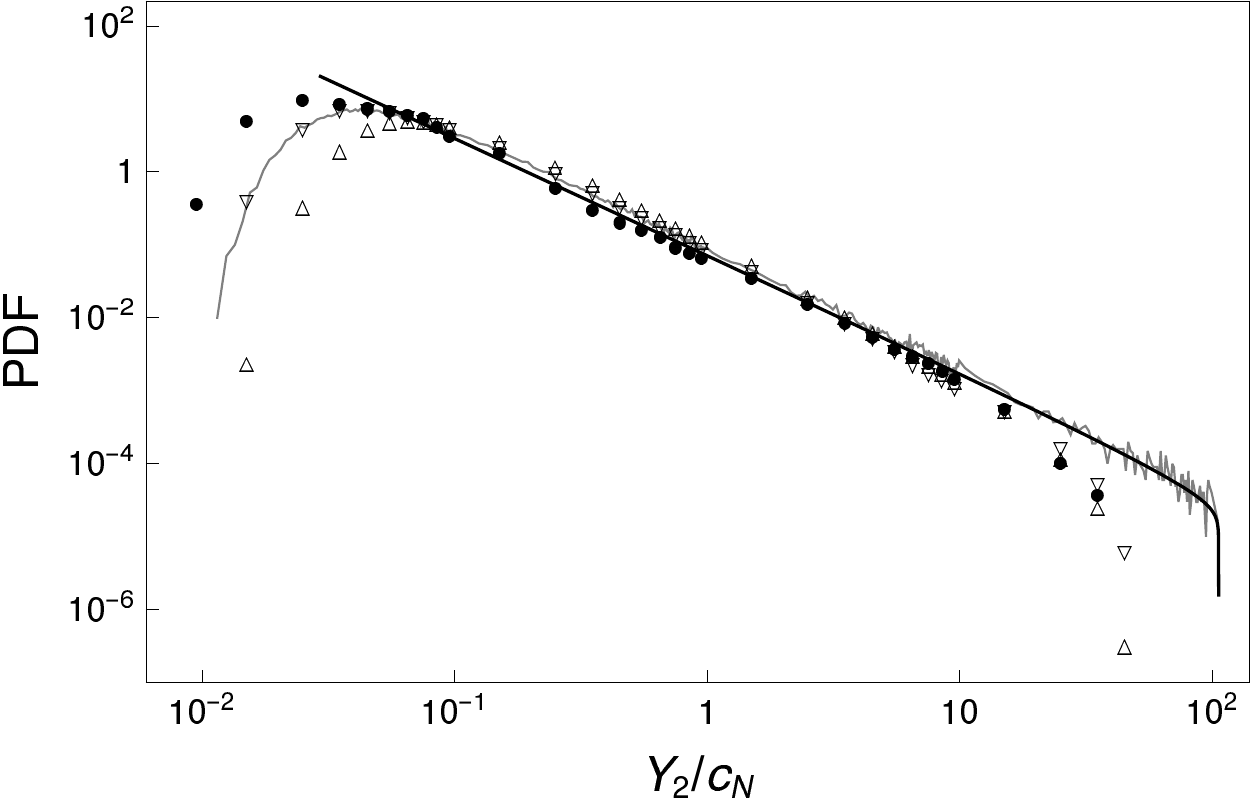}&
      \includegraphics[height=.3\textwidth,viewport=0 0 360 233]{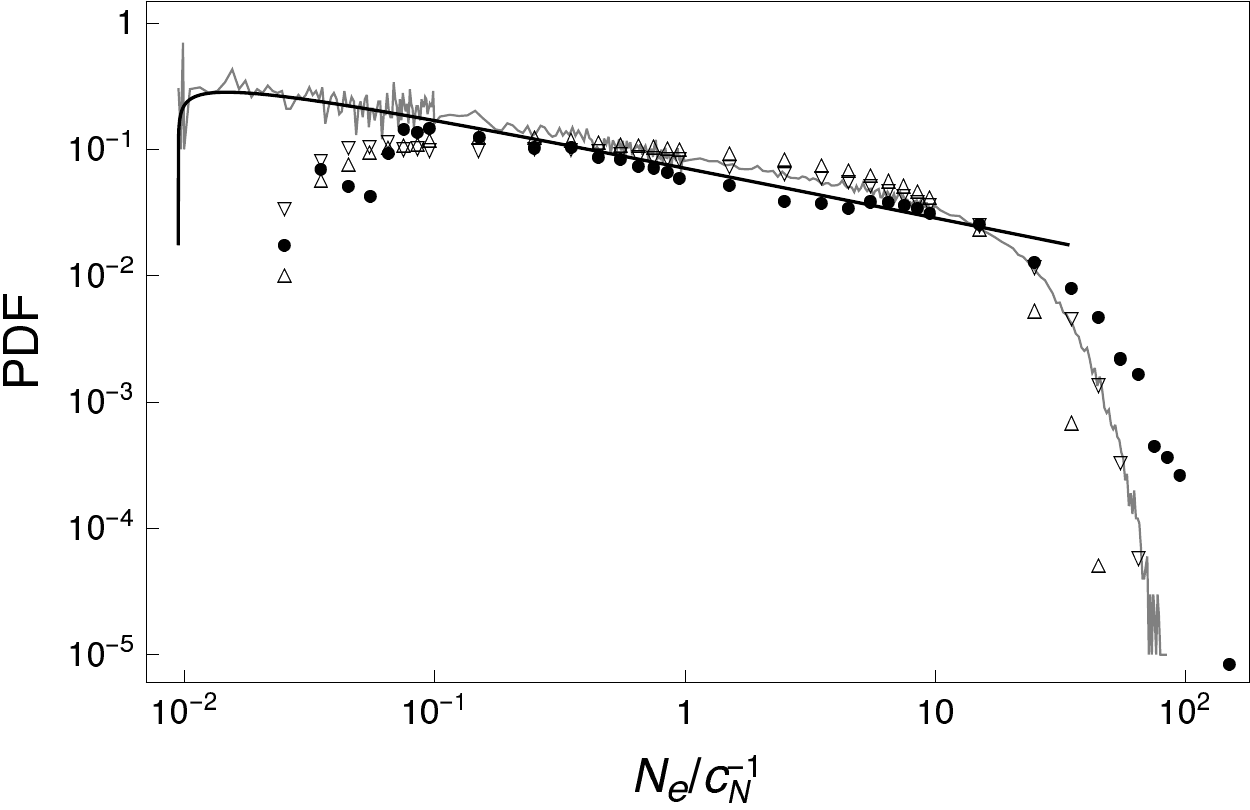}
 \end{tabular}
 \caption{\small
Probability distributions.
(a)~$\mathrm{\Pi}_{W_{1,N}}$ of the largest weight,
(b)~$\mathrm{\Pi}_{W_{2,N}}$ of the second largest weight,
(c)~$\mathrm{\Pi}_{Y_2}$ of the sample mean weight of families,
and (d)~$\mathrm{\Pi}_{N_{\mathrm{e}}}$ of the effective population size,
obtained from the recursions (Eqs.~\ref{eqn:recursion-Y2}, \ref{eqn:recursion-W1} and~\ref{eqn:recursion-W2})
with $\alpha=1.2$ and
$N=10^4{\;}(\vartriangle)$, $10^5{\;}(\triangledown)$, and $10^6{\;}(\bullet)$.
These are compared to histograms (depicted by noisy gray curves), obtained by Pareto sampling ($N=10^6$).
The black lines are from Eqs.\eqref{eqn:largest-W}, \eqref{eqn:2nd-largest-W}, \eqref{eqn:distr-Y2}, and~\eqref{eqn:distr-Ne} with $N=10^6$.
$Y_2$ and $N_{\mathrm{e}}$ are scaled by $c_{\scalebox{0.55}{$N$}}$ and $c_{\scalebox{0.55}{$N$}}^{-1}$, respectively.
 }\label{fig:derrida-trick}
\end{figure}

I use a trick given in \cite{Derrida-Flyvbjerg87,Derrida-Flyvbjerg87-PhysA} to generate distributions of $W_{1,N}$, $W_{2,N}$, and $Y_2$ for $1<\alpha<2$.
Suppose one has a set $\qty{W_1,\ldots,W_{N-1}}$ of $(N-1)$ weights of families and a corresponding value of $Y_2$.
Let $W_{\mathrm{max}}$ and $W'_{\mathrm{max}}$ be the two largest values of this set of weights.
One can add another family of weight $W_{N}$ with probability distribution $N^{-1}\rho(w)$ where $\rho(w)$ is as in Eq.\eqref{eqn:freq-spectrum-beta}, and simultaneously shrink all the other weights by a factor $(1-W_N)$.
One now has a new set with a new value given by the recursion relation
\begin{equation}\label{eqn:Derrida-tric}
 \tilde{Y}_2=W_N^2+(1-W_N)^2 Y_2.
\end{equation}
If $\tilde{W}_{\mathrm{max}}$ and $\tilde{W}'_{\mathrm{max}}$ are the two largest weights of the new set, then one has the recursions
\begin{align*}
 \tilde{W}_{\mathrm{max}}&=\max\qty(W_N,(1-W_N)W_{\mathrm{max}})\\
 \tilde{W}'_{\mathrm{max}}&=\max\qty(\min\qty(W_N,(1-W_N)W_{\mathrm{max}}),(1-W_N)W'_{\mathrm{max}}).
\end{align*}

I just calculate a random sequence $\qty{y_k}$ ($k=0,1,2,\ldots$) with $y_0$ being randomly chosen between 0 and 1, where the recursion relation which gives the $y_k$'s is
\begin{equation}\label{eqn:recursion-Y2}
 y_{k+1}=W_k^2+(1-W_k)^2 y_k
\end{equation}
and where the $W_k$'s are randomly chosen between $\varepsilon_{\scalebox{0.55}{$N$}}$ and 1 according to $\rho(w)$, where $\varepsilon_{\scalebox{0.55}{$N$}}$ is as in Eq.\eqref{eqn:cut-off}.
The $y_k$'s (resp. $y_k^{-1}$'s), generated by this iterative procedure,
are distributed according to $\mathrm{\Pi}_{Y_2}$ (resp. $\mathrm{\Pi}_{N_{\mathrm{e}}}$).
One can also construct in the same way sequences of $w_k$ and $w'_k$ by iterative procedures
\begin{align}
 w_{k+1}&=\max\qty(W_k,(1-W_k)w_k)
 \label{eqn:recursion-W1}\\
 w'_{k+1}&=\max\qty(\min\qty(W_k,(1-W_k)w_k),(1-W_k)w'_k)
 \label{eqn:recursion-W2}
\end{align}
where again the $W_k$'s are distributed according to $\rho(w)$.
By drawing the histograms of $w_k$'s, $w'_k$'s, $y_k$'s and $y_k^{-1}$'s,
one can get
the probabilities $\mathrm{\Pi}_{W_{1,N}}$, $\mathrm{\Pi}_{W_{2,N}}$, $\mathrm{\Pi}_{Y_2}$ and $\mathrm{\Pi}_{N_{\mathrm{e}}}$, respectively.
Fig.~\ref{fig:derrida-trick} shows the results obtained after $2\times 10^8$ iterations with parameters being set to $\alpha=1.2$, and with populations of $N=10^4$, $10^5$ and $10^6$.

The simulation results of the Pareto sampling are also shown in Fig.~\ref{fig:derrida-trick},
where $X_i$'s ($N=10^6$ trials) are drawn from the Pareto distribution (Eq.~\ref{eqn:Pareto-density}) with $\alpha=1.2$,
and the histograms of $W_{1,N}$'s, $W_{2,N}$'s, $Y_2$'s and $Y_2^{-1}$'s
are obtained from $10^5$ independent realizations.

\section{Genetic variation in marine populations}\label{sect:application}
Highly fecund marine species undergo large and intermittent fluctuations in recruitment
\cite{Hjort1914}.
Very large interfamilial, or sweepstakes, variation in reproductive success has been documented
\cite{Hedgecock-Pudovkin2011}.
While sweepstakes reproduction appears to be prevalent in marine systems,
it is somewhat paradoxical to have observed the absence of reduced genetic diversity from a single reproduction event
\cite{Ruzzante96,Flowers2002,Rose-etal2006,Selkoe-etal2006,Buston-etal2009,Christie-etal2010,Hogan-meps2010,He-atal2012,Jue2014,Cornwell2016,Riquet2017}.
Even a species that showed evidence for sweepstakes in one place may not show evidence for it at another time or in another place
\cite{Hedgecock-etal2007,Taris-etal2009}.
This limited evidence can be attributed to variation in intensity of sweepstakes reproduction across time, that is, ``the right place, but the wrong time'' as stated in \cite{Hedgecock-Pudovkin2011},
but it has not yet been fully proved or explained.

\begin{figure}[hbt!]
 \centering
 \includegraphics[height=.3\textwidth,viewport=0 0 360 231]{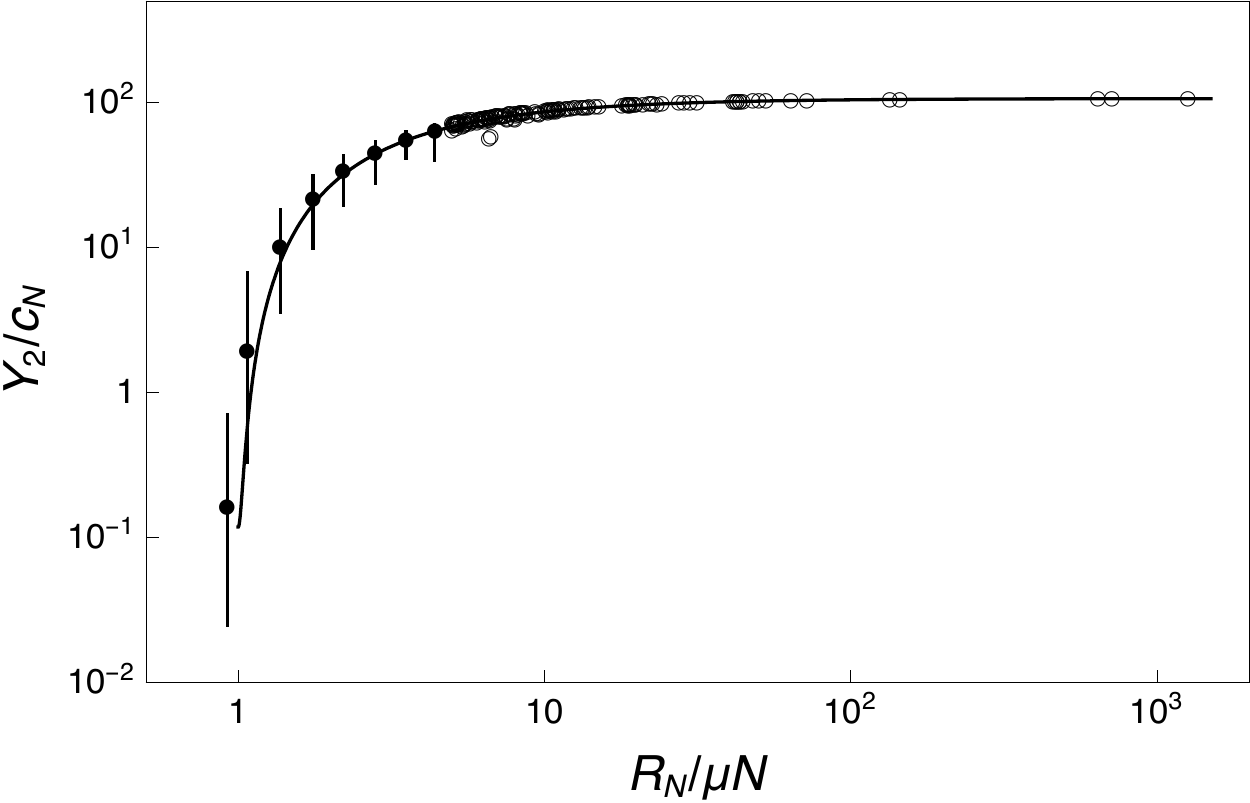}
  \caption{\small
The more recruitment, the more genetic drift.
The $Y_2$'s are shown associated with recruitment strength $R_N$,
after running $10^5$ repetitions of random Pareto sampling ($\alpha=1.2$ and $N=10^6$).
Solid circles and error bars mark the mean and 95\% of binned data.
The solid curve is from Eq.\eqref{eqn:sweepstakes}.
 $Y_2$ and $R_N$ are scaled by $c_{\scalebox{0.55}{$N$}}$ and $\mu N$, respectively.
 }\label{diversity-loss}
\end{figure}

As an application example, the concept of the statistical domination by the largest term is used to assess an increasing reduction in genetic variation at larger recruitment $R_N$.
$(N-1)$ lower order statistics, $X_{2,N}\geq\cdots\geq X_{N,N}$, of $N$ independent $\mathrm{Pareto}(\alpha)$ random variables with $1<\alpha<2$ have finite first and second moments.
Letting $U_2$ be the sum of squared weights of these $(N-1)$ lower order statistics,
\begin{equation*}
 U_2=\sum_{i=2}^N\qty({X_{i,N}}/{R_{2,N}})^2,
\end{equation*}
one has, from Eq.\eqref{eqn:Derrida-tric},
\begin{equation}\label{eqn:sweepstakes}
 \begin{aligned}
 Y_2 &= W_{1,N}^2 + (1-W_{1,N})^2 U_2\\
 &= \qty(1-\frac{\mu N}{R_N})^2+
 \qty(\frac{\mu N}{R_N})^2
 \frac{(\alpha-1)^2\mathrm{\Gamma}(2-2/\alpha)N^{2(1/\alpha-1)}}{\alpha(\alpha-2)}
 \end{aligned}
\end{equation}
with the second term being negligible.
Eq.\eqref{eqn:sweepstakes} agrees with the result from the simulation of random Pareto sampling (performed in \S\ref{section:Derrida-trick}), as shown in Fig.~\ref{diversity-loss}.
The typical value of $Y_2/c_{\scalebox{0.55}{$N$}}$ is of order
$N^{(\alpha-1)(\alpha-2)/\alpha}$
($=10^{-0.8}$ with $\alpha=1.2$ and $N=10^6$).
One sees the non-self-averaging effect, that is, $Y_2$ depends on the recruit sample.
A large reduction in genetic variation will occur in years of large recruitment.

\section{Conclusions}
When tracing a history of realizations of family-size frequencies in a population with $\mathrm{Pareto}(\alpha)$ offspring-number distribution of index $1<\alpha<2$,
one may take the average of the weights of families over the whole population existing at any given time.
The statistical domination by the largest family leads to large fluctuations in the population average.
So the non-self-averaging behavior emerges.

I have studied the fluctuations for the largest and second largest weights of families,
and the fluctuations for the average weight $Y_2$ and for its reciprocal (i.e. the effective population size $N_{\mathrm{e}}$).
I have obtained asymptotic expressions for the probability distributions of these quantities.
The most probable value of the $Y_2$ is close to the typical value of squared weight of the largest family ($\sim N^{2/\alpha-2}$),
and differs from its mean value $\mathrm{E}[Y_2]{\,}(\sim N^{1-\alpha})$ of the reproduction process.
The $N_{\mathrm{e}}$ has a broad distribution, with an upper cut-off corresponding to the inverse square of the typical largest weight.
The most probable $N_{\mathrm{e}}=(2-\alpha)^{-2}$, independent of the population size $N$,
is close to the lower bound of the distribution
and very different from the harmonic mean of $N_{\mathrm{e}}$'s over replicate populations.
There is a broken symmetry for scaling of the typical $Y_2$ and of its reciprocal.

Non-self-averaging effects are crucial in understanding the complexities surrounding intermittent, large recruitment events typical in marine populations.
In occasional years of large recruitment,
only a few parents will contribute a dominant fraction of the recruitment.
Random changes in genetic diversity are associated with recruitment strength, and large recruitment events enhance genetic drift with reduced $N_{\mathrm{e}}$.

\bibliographystyle{unsrt}
\bibliography{bib-niwa-genetics}

\newcommand{\SortNoop}[1]{}
\begin{thebibliography}{10}

\bibitem{SBB13}
M.~Steinr{\"u}cken, M.~Birkner, and J.~Blath.
\newblock Analysis of {DNA} sequence variation within marine species using
  {Beta}-coalescents.
\newblock {\em Theor. Popul. Biol.}, 87:15--24, 2013.

\bibitem{Arnason-Halldorsdottir2015}
E.~{\'A}rnason and K.~Halld{\'o}rsd{\'o}ttir.
\newblock Nucleotide variation and balancing selection at the \textit{{Ckma}}
  gene in {Atlantic} cod: analysis with multiple merger coalescent models.
\newblock {\em PeerJ}, 3:e786, 2015.

\bibitem{Niwa-etal2016}
H.-S. Niwa, K.~Nashida, and T.~Yanagimoto.
\newblock Reproductive skew in {Japanese} sardine inferred from {DNA}
  sequences.
\newblock {\em ICES J. Mar. Sci.}, 73:2181--2189, 2016.

\bibitem{Reed-Hughes2002}
W.~J. Reed and B.~D. Hughes.
\newblock From gene families and genera to incomes and internet file sizes: why
  power laws are so common in nature.
\newblock {\em Phys. Rev. E}, 66:067103, 2002.

\bibitem{Niwa-etal2017}
H.-S. Niwa, K.~Nashida, and T.~Yanagimoto.
\newblock Allelic inflation in depleted fish populations with low recruitment.
\newblock {\em ICES J. Mar. Sci.}, 74:1639--1647, 2017.

\bibitem{Derrida-PhysD97}
B.~Derrida.
\newblock From random walks to spin glasses.
\newblock {\em Physica D}, 107:186--198, 1997.

\bibitem{Wright1931}
S.~Wright.
\newblock Evolution in {Mendelian} populations.
\newblock {\em Genetics}, 16:96--159, 1931.

\bibitem{Schweinsberg2003}
J.~Schweinsberg.
\newblock Coalescent processes obtained from supercritical {Galton-Watson}
  processes.
\newblock {\em Stoch. Process. Their Appl.}, 106:107--139, 2003.

\bibitem{Huillet2014}
T.~Huillet.
\newblock Pareto genealogies arising from a {Poisson} branching evolution model
  with selection.
\newblock {\em J. Math. Biol.}, 68:727--761, 2014.

\bibitem{MezardPSTV84}
M.~M\'ezard, G.~Parisi, N.~Sourlas, G.~Toulouse, and M.~Virasoro.
\newblock Replica symmetry breaking and the nature of the spin glass phase.
\newblock {\em J. Physique}, 45:843--854, 1984.

\bibitem{Derrida-Peliti91}
B.~Derrida and L.~Peliti.
\newblock Evolution in a flat fitness landscape.
\newblock {\em Bull. Math. Biol.}, 53:355--382, 1991.

\bibitem{Serva2004}
M.~Serva.
\newblock Lack of self-averaging and family trees.
\newblock {\em Physica A}, 332:387--393, 2004.

\bibitem{Serva2005}
M.~Serva.
\newblock On the genealogy of populations: trees, branches and offspring.
\newblock {\em J. Stat. Mech.}, 2005:P07011, 2005.

\bibitem{Bolthausen-Sznitman98}
E.~Bolthausen and A.-S. Sznitman.
\newblock On {Ruelle's} probability cascades and an abstract cavity method.
\newblock {\em Commun. Math. Phys.}, 197:247--276, 1998.

\bibitem{Brunet-Derrida2009}
{\'{E}}.~Brunet and B.~Derrida.
\newblock Statistics at the tip of a branching random walk and the delay of
  traveling waves.
\newblock {\em Europhys. Lett.}, 87:60010, 2009.

\bibitem{Brunet-Derrida2011}
{\'{E}}.~Brunet and B.~Derrida.
\newblock A branching random walk seen from the tip.
\newblock {\em J. Stat. Phys.}, 143:420--446, 2011.

\bibitem{Brunet-Derrida2013}
{\'{E}}.~Brunet and B.~Derrida.
\newblock Genealogies in simple models of evolution.
\newblock {\em J. Stat. Mech.}, 2013:P01006, 2013.

\bibitem{Stewart76}
F.~M. Stewart.
\newblock Variability in the amount of heterozygosity maintained by neutral
  mutations.
\newblock {\em Theor. Popul. Biol.}, 9:188--201, 1976.

\bibitem{Higgs-PRE95}
P.~G. Higgs.
\newblock Frequency distributions in population genetics parallel those in
  statistical physics.
\newblock {\em Phys. Rev. E}, 51:95--101, 1995.

\bibitem{Tajima83}
F.~Tajima.
\newblock Evolutionary relationship of {DNA} sequences in finite populations.
\newblock {\em Genetics}, 105:437--460, 1983.

\bibitem{King2018}
L.~King, J.~Wakeley, and S.~Carmi.
\newblock A non-zero variance of {Tajima's} estimator for two sequences even
  for infinitely many unlinked loci.
\newblock {\em Theor. Popul. Biol.}, 122:22--29, 2018.

\bibitem{Fisher-etal43}
R.~A. Fisher, A.~Steven Corbet, and C.~B. Williams.
\newblock The relation between the number of species and the number of
  individuals in a random sample of an animal population.
\newblock {\em J. Anim. Ecol.}, 12:42--58, 1943.

\bibitem{Watterson74}
G.~A. Watterson.
\newblock Models for the logarithmic species abundance distributions.
\newblock {\em Theor. Popul. Biol.}, 6:217--250, 1974.

\bibitem{Levy1937}
P.~L{\'e}vy.
\newblock {\em Th{\'e}eorie de l'Addition des Variables Al{\'e}atoires}.
\newblock Gauthier-Villars, Paris, 1937.

\bibitem{ZKS05}
I.~V. Zaliapin, Y.~Y. Kagan, and F.~P. Schoenberg.
\newblock Approximating the distribution of {Pareto} sums.
\newblock {\em Pure Appl. Geophys.}, 162:1187--1228, 2005.

\bibitem{Bouchaud-Georges90}
J.~P. Bouchaud and A.~Georges.
\newblock Anomalous diffusion in disordered media: statistical mechanisms,
  models and physical applications.
\newblock {\em Phys. Rep.}, 195:127--293, 1990.

\bibitem{Hofstad2016}
R.~van~der Hofstad.
\newblock {\em Random Graphs and Complex Networks}.
\newblock Cambridge University Press, Cambridge, UK., 2016.

\bibitem{Hausdorff1923}
F.~Hausdorff.
\newblock Momentprobleme f{\"u}r ein endliches {Intervall}.
\newblock {\em Math. Z.}, 16:220--248, 1923.

\bibitem{Derrida-Toulouse85}
B.~Derrida and G.~Toulouse.
\newblock Sample to sample fluctuations in the random energy model.
\newblock {\em J. Physique Lett.}, 46:223--228, 1985.

\bibitem{Pitman-Yor97}
J.~Pitman and M.~Yor.
\newblock The two-parameter poisson-dirichlet distribution derived from a
  stable subordinator.
\newblock {\em Ann. Probab.}, 25:855--900, 1997.

\bibitem{Derrida-Flyvbjerg87}
B.~Derrida and H.~Flyvbjerg.
\newblock The random map model: a disordered model with deterministic dynamics.
\newblock {\em J. Physique}, 48:971--978, 1987.

\bibitem{Derrida-Flyvbjerg87-PhysA}
B.~Derrida and H.~Flyvbjerg.
\newblock Statistical properties of randomly broken objects and of multivalley
  structures in disordered systems.
\newblock {\em J. Phys. A: Math. Gen.}, 20:5273--5288, 1987.

\bibitem{Romeo2003}
M.~Romeo, V.~Da~Costa, and F.~Bardou.
\newblock Broad distribution effects in sums of lognormal random variables.
\newblock {\em Eur. Phys. J. B}, 32:513--525, 2003.

\bibitem{Hjort1914}
J.~Hjort.
\newblock Fluctuations in the great fisheries of the northern {Europe} viewed
  in the light of biological research.
\newblock {\em Rapp. P.-V. R{\'e}un. Cons. Int. Explor. Mer}, 20:1--228, 1914.

\bibitem{Hedgecock-Pudovkin2011}
D.~Hedgecock and A.~I. Pudovkin.
\newblock Sweepstakes reproductive success in highly fecund marine fish and
  shellfish: a review and commentary.
\newblock {\em Bull. Mar. Sci.}, 87:971--1002, 2011.

\bibitem{Ruzzante96}
D.~E. Ruzzante, C.~T. Taggart, and D.~Cook.
\newblock Spatial and temporal variation in the genetic composition of a larval
  cod (\textit{Gadus morhua}) aggregation: cohort contribution and genetic
  stability.
\newblock {\em Can. J. Fish. Aquat. Sci.}, 53:2695--2705, 1996.

\bibitem{Flowers2002}
J.~M. Flowers, S.~C. Schroeter, and R.~S Burton.
\newblock The recruitment sweepstakes has many winners: genetic evidence from
  the sea urchin \textit{Strongylocentrotus purpuratus}.
\newblock {\em Evolution}, 56:1445--1453, 2002.

\bibitem{Rose-etal2006}
C.~G. Rose, K.~T. Paynter, and M.~P. Hare.
\newblock Isolation by distance in the eastern oyster, \textit{Crassostrea
  virginica}, in {Chesapeake Bay}.
\newblock {\em J. Hered.}, 97:158--170, 2006.

\bibitem{Selkoe-etal2006}
K.~A. Selkoe, S.~D. Gaines, J.~E. Caselle, and R.~R. Warner.
\newblock Current shifts and kin aggregation explain genetic patchiness in fish
  recruits.
\newblock {\em Ecology}, 87:3082--3094, 2006.

\bibitem{Buston-etal2009}
P.~M. Buston, C.~Fauvelot, M.~Y.~L. Wong, and S.~Planes.
\newblock Genetic relatedness in groups of the humbug damselfish
  \textit{Dascyllus aruanus}: small, similar-sized individuals may be close
  kin.
\newblock {\em Mol. Ecol.}, 18:4707--4715, 2009.

\bibitem{Christie-etal2010}
M.~R. Christie, D.~W. Johnson, C.~D. Stallings, and M.~A. Hixon.
\newblock Self-recruitment and sweepstakes reproduction amid extensive gene
  flow in a coral-reef fish.
\newblock {\em Mol. Ecol.}, 19:1042--1057, 2010.

\bibitem{Hogan-meps2010}
J.~D. Hogan, R.~J. Thiessen, and D.~D. Heath.
\newblock Variability in connectivity indicated by chaotic genetic patchiness
  within and among populations of a marine fish.
\newblock {\em Mar. Ecol. Prog. Ser.}, 417:263--275, 2010.

\bibitem{He-atal2012}
Y.~He, S.~E. Ford, D.~Bushek, E.~N. Powell, Z.~Bao, and X.~Guo.
\newblock Effective population sizes of eastern oyster \textit{Crassostrea
  virginica} ({Gmelin}) populations in {Delaware Bay, USA}.
\newblock {\em J. Mar. Res.}, 70:357--379, 2012.

\bibitem{Jue2014}
N.~K. Jue, F.~C. Coleman, and C.~C. Koenig.
\newblock Wide-spread genetic variability and the paradox of effective
  population size in the gag, \textit{Mycteroperca microlepis}, along the {West
  Florida Shelf}.
\newblock {\em Mar. Biol.}, 161:1905--1918, 2014.

\bibitem{Cornwell2016}
B.~H. Cornwell, J.~L. Fisher, S.~G. Morgan, and J.~E. Neigel.
\newblock Chaotic genetic patchiness without sweepstakes reproduction in the
  shore crab \textit{Hemigrapsus oregonensis}.
\newblock {\em Mar. Ecol. Prog. Ser.}, 548:139--152, 2016.

\bibitem{Riquet2017}
F.~Riquet, T.~Comtet, T.~Broquet, and F.~Viard.
\newblock Unexpected collective larval dispersal but little support for
  sweepstakes reproductive success in the highly dispersive brooding mollusc
  \textit{Crepidula fornicata}.
\newblock {\em Mol. Ecol.}, 26:5467--5483, 2017.

\bibitem{Hedgecock-etal2007}
D.~Hedgecock, S.~Launey, A.~I. Pudovkin, Y.~Naciri, S.~Lap{\`e}gue, and
  F.~Bonhomme.
\newblock Small effective number of parents ({$N_b$}) inferred for a naturally
  spawned cohort of juvenile {European} flat oysters \textit{Ostrea edulis}.
\newblock {\em Mar. Biol.}, 150:1173--1182, 2007.

\bibitem{Taris-etal2009}
N.~Taris, P.~Boudry, F.~Bonhomme, M.~D. Camara, and S~Lapegue.
\newblock Mitochondrial and nuclear {DNA} analysis of genetic heterogeneity
  among recruitment cohorts of the {European} flat oyster \textit{Ostrea
  edulis}.
\newblock {\em Biol. Bull.}, 217:233--241, 2009.

\end{thebibliography}

\end{document}